\documentclass[12pt]{amsart}
\usepackage{amsmath,amssymb, euscript, graphicx, hyperref}

\newtheorem{thm}{Theorem}[section]
\newtheorem{lem}[thm]{Lemma}
\newtheorem{prop}[thm]{Proposition}
\theoremstyle{definition}
\newtheorem{defn}[thm]{Definition}

\newtheorem{rem}[thm]{Remark}
\newtheorem{cor}[thm]{Corollary}

\newcommand{\blackboard}[1]{\ensuremath{\mathbb{#1}}}

\newcommand{\Z}{\blackboard{Z}}

\begin{document}

\address{Azer Akhmedov, Department of Mathematics,
North Dakota State University,
Fargo, ND, 58102, USA}
\email{azer.akhmedov@ndsu.edu}

\address{Cody Martin, Department of Mathematics,
North Dakota State University,
Fargo, ND, 58102, USA}
\email{cody.martin@ndsu.edu}

 \begin{center} {\bf {\Large Non-bi-orderability of $6_2$ and $7_6$}} \end{center}
 
 \bigskip
 
  \begin{center} {\bf Azer Akhmedov and Cody Martin} \end{center}
 
 \bigskip
 
 {\small ABSTRACT: We prove that the knot groups of $6_2$ and $7_6$ are not bi-orderable. These are the only two knot groups up to 7 crossings whose bi-orderability was not known. Our method applies to a very broad class of knots.\footnote{In our next publication, we will present a complete characterization of the bi-orderability of knot groups. This characterization (a necessary and sufficient condition) does not produce a finite algorithm for deciding the bi-orderabilty of an arbitrary knot group but at the practical level (possibly with the help of computers) it is effective in treating a very broad class of knots.}

 \section{Introduction}

 Given a knot $K \subset \mathbb{S}^3$, the group $\pi _1(\mathbb{S}^3 \backslash K)$ is called the {\em knot group} of $K$. It is well known that every knot group is left-orderable \cite {HS, CR}. On the other hand, not every knot group is bi-orderable; for example, it is not difficult to see that the knot group of a non-trivial torus knot is not bi-orderable. Indeed, the fundamental group of the $(p,q)$-torus knot complement has a presentation $\langle x, y \ | \ x^p = y^q \rangle $ where $p, q\geq 2$. Thus a power of $x$ commutes with $y$, while it is easy to show that $x$ does not commute with $y$. However, this is impossible in a bi-orderable group: assume that $<$ is a bi-order in such a group. Then either $xyx^{-1} > y$ or $xyx^{-1} < y$. Without loss of generality assume that $xyx^{-1} > y$. Then, inductively, we obtain $x^nyx^{-n} > x^{n-1}yx^{1-n}$ for every natural number $n$. Taking $n = p$, we obtain a contradiction.
 
 \medskip
 
 The well known argument from the previous paragraph shows that the groups $\pi _1(\mathrm{{\bf Klein \ Bottle}})$ and $\pi _1(\mathrm{{\bf Trefoil}})$ are not bi-orderable. On the other hand, the Figure-8 knot  turns out to be bi-orderable. It is possible, although not easy, to see this directly from the following presentation of the Figure-8 knot group: $$\langle t, a, b \ | \ tat^{-1} = aba, tbt^{-1} = ab \rangle $$
 
 \medskip
 
 In \cite {PR} a very interesting criterion (a sufficient condition) has been proved for the bi-orderability of knots. 
 
 \begin{thm}\label{thm:CR} Let $K$ be a fibered knot such that all the roots of the Alexander polynomial $\Delta _K(t)$ are positive and real. Then $\pi _1(K)$ is bi-orderable. 
 \end{thm}
 
 \medskip
 
 The unknot, the trefoil knot, and the Figure-8 knot are all fibered. The Alexander polynomial of the Figure-8 knot is given by $\Delta (t) = t^2-3t+1$ thus all the roots are indeed real and positive. For the Trefoil knot, however, the Alexander polynomial equals $\Delta (t) = t^2-t+1$ which has imaginary roots.
 
 \medskip
 
 A major ingredient in the proof of Theorem \ref{thm:CR} is the involvement of Seifert-van Kampen Theorem by using the fiberedness of the knot. By far, not all knots are fibered (although both of the knots $6_2$ and $7_6$ that we treat in this paper are indeed fibered). In \cite{CR1} a partial converse to Theorem \ref{thm:CR} is given which states that if a fibered knot group is bi-orderable then $\Delta _K(t)$ has at least one positive real root. In \cite{CR}, using the results of Chiswell-Glass-Wilson \cite{CGW}, the bi-orderability and non-bi-orderability of several non-fibered knots has also been deduced. It has also been asked in \cite{CR1} if the full converse of their main result holds, i.e. if the Alexander polynomial of a fibred knot $K$ has at least one real positive root then is it true that $K$ is bi-orderable? {\em The result of our paper settles this question in the negative.}
 
 \medskip
 
  In this paper, we involve a Seifert surface technique to obtain nice presentations for the knots $6_2$ and $7_6$. We borrow this technique from \cite{Li} and \cite{FSW}; it applies to any knot, and allows to write a knot group as an HNN extension (of a free group). Moreover, it allows to drop the fiberdness condition on the knots ({\em this indeed settles another question asked in \cite{CR1})}. Then we use the concept of subgroups of infinitesimals; we develop this concept (which is interesting independently) to make it applicable to the knot groups in question. Our method applies to a very broad class of knots. In the current paper, we treat only the knots $6_2$ and $7_6$ which were the only two knots up to crossing number 7 whose bi-orderability was not known.   
 
 \medskip
 
 Given any knot $K \subset \mathbb{S}^3$ one can find an embedded surface $\Sigma \subset \mathbb{S}^3$ which bounds $K$. For example, both a disk $\mathbb{D}^2$ and a M\"obius Band bounds an unknot, so does any punctured compact surface. In general, there are many non-homeomorphic surfaces bounding a given knot. One can indeed find such surfaces with some nice properties. For example, it is always true that there exists an {\em orientable surface} $\Sigma $ which bounds $K$. Such a surface $\Sigma $ is called  {\em a Seifert surface of the knot $K$}. The existence of such surfaces was first proved by Pontriagin and Frankl in 1930. In 1934, Seifert gave another proof also providing a rather nice and short algorithm for constructing such a surface. 
 
 \begin{figure}[h!]
  \includegraphics[width=6in,height=5in]{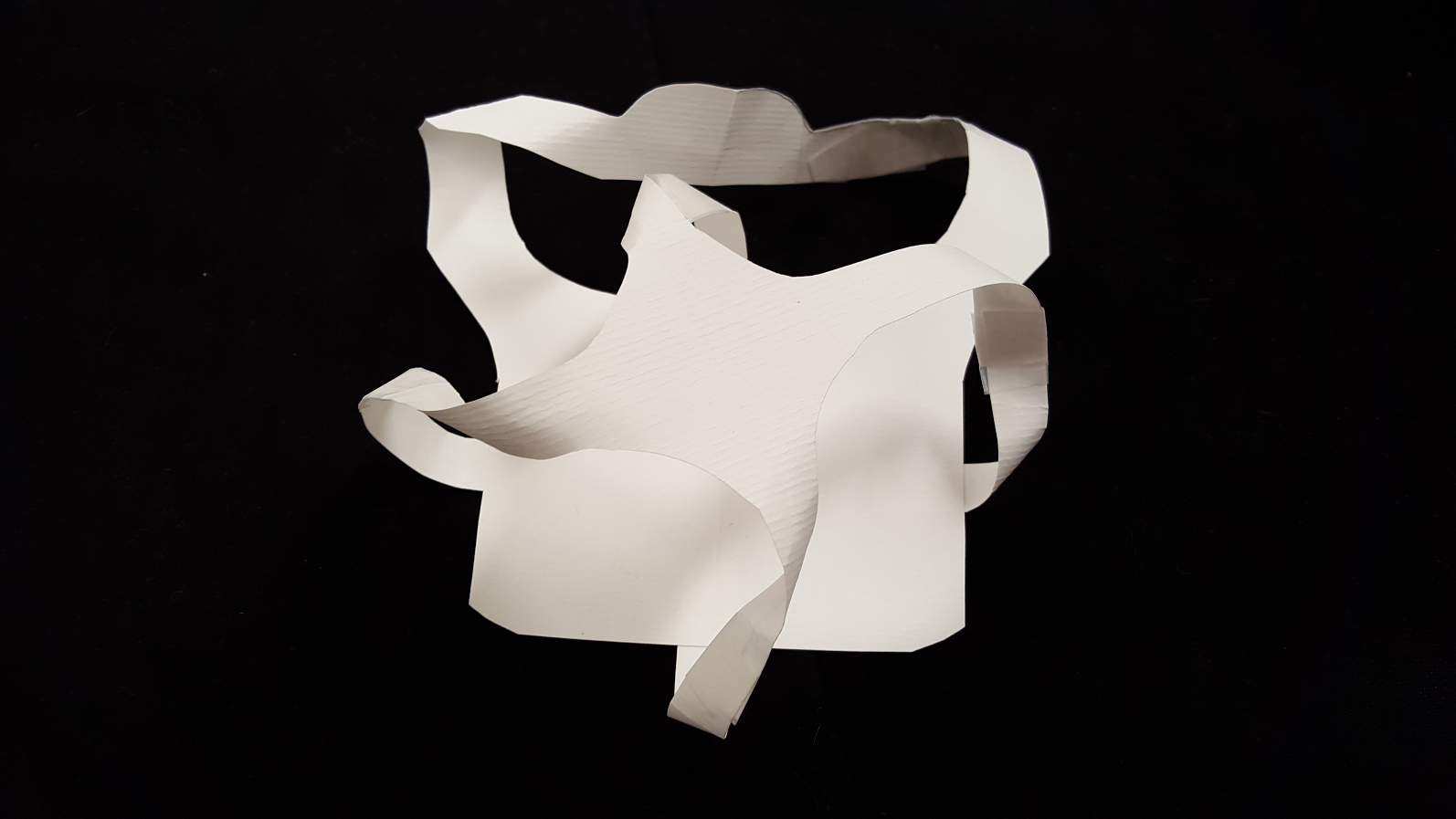}
\caption{A Seifert surface of $6_2$.}
\label{labelname}
\end{figure}
 
 \medskip
 
 Seifert surfaces are never unique since one can always increase the genus of it without changing the boundary. Within a fixed homeomorphism type, there is still a rich variety of possibilities for Seifert surfaces; for example, the isomorphism type of the group $\pi _1(\mathbb{S}^3 \backslash \Sigma )$ may vary depending on the surface. If this group is free then $\Sigma $ is called {\em a free Seifert surface} or {\em an unknotted Seifert surface}, otherwise it is called {\em a non-free Seifert surface} or {\em a knotted Seifert surface}. A Seifert surface may or may not be incompressible. Every knot admits an incompressible Seifert surface; by Loop Theorem, any Seifert surface with minimal genus is necessarily incompressible. Also, every knot admits both free and non-free Seifert surfaces if one allows compressible surfaces \cite{O}. On the other hand, there exist knots which admit only non-free incompressible Seifert surfaces \cite{L}. Let us also mention that the Seifert surface obtained as a result of the Seifert algorithm is always free.
 
 \medskip
 
 Let $K\subset \mathbb{S}^3$ be a knot, $\nu (K)$ be its tubular neighborhood. Following [FSW], let $\Sigma $ be an incompressible Seifert surface of $K$ and $X(K) = \mathbb{S}^3\backslash \nu (K)$. By taking a tubular neighborhood $\Sigma \times (-1, 1)$ we view $X(K) \backslash (\Sigma \times (-1, 1))$ as a manifold obtain by cutting along $\Sigma $. One can glue the positive and negative copies $\Sigma \times \{1\}$ and $\Sigma \times \{-1\}$ to obtain the original manifold $X(K)$, thus Seifert-van Kampen Theorem gives us the following presentation: $$\pi _1(X(K)) = \langle \pi _1(X(K)\backslash (\Sigma \times (-1, 1))), t \ | \ tx_+t^{-1} = x_-, x_+\in \pi _1(\Sigma \times \{1\})\rangle $$ where $x_-\in \pi _1(\Sigma \times \{-1\})$ is the parallel copy of $x_+$ ``on the other side" of $\Sigma $.  
 
 \medskip
 
 The above presentation often turns out to be quite pleasant for computational purposes (even though it takes some work to compute this presentation), especially when compared with the Wirtinger presentation of the knot complement. The presentation of the Figure-8 knot group mentioned above can be obtained by this technique. Applying it to the Trefoil knot, one obtains the following well known presentation: $$\langle t, a, b \ | \ tat^{-1} = b, tbt^{-1} = ba^{-1} \rangle $$

 We will obtain such presentations for the fundamental groups of the knot complements of $6_2$ and $7_6$. Both of these knot groups will have a structure of a semidirect product $\mathbb{Z}\ltimes \mathbb{F}_4$.          
    
 \begin{figure}[h!]
  \includegraphics[width=6in,height=5in]{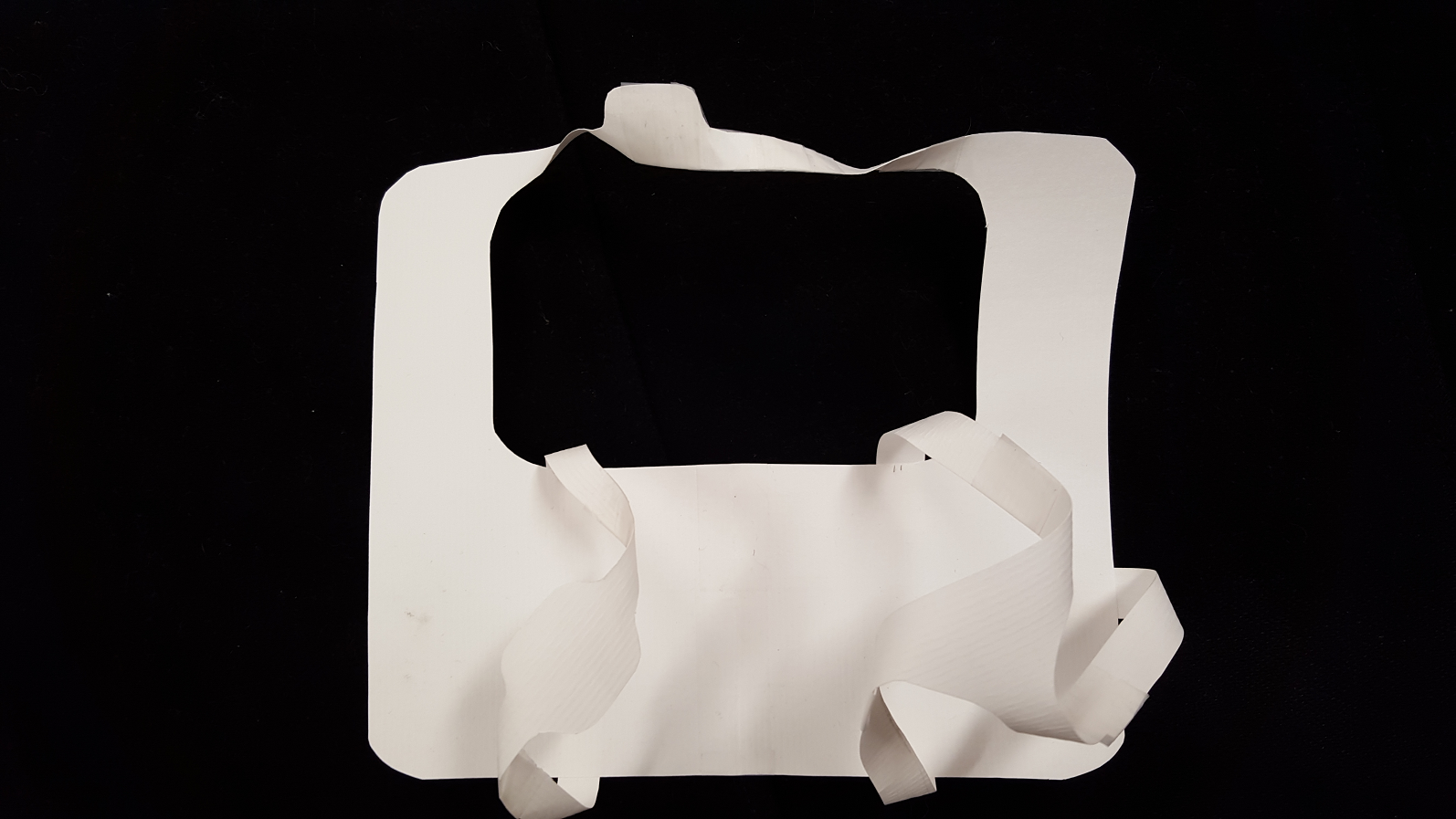}
\caption{A Seifert surface of $7_6$.}
\label{labelname}
\end{figure}

\bigskip
 
 \section{Subgroup of Infinitesimals}

 In this section we introduce a very useful notion of a {\em subgroup of infinitesimals} and study its properties. 
 
 \medskip
 
 Let $\Gamma $ be a group and $<$ be a left-order in it. For any $x\in \Gamma $, we let  $|x| = x$ if $x > 1$ and $|x| = x^{-1}$ otherwise. For $f, g\in \Gamma $ we will say $f$ is infinitesimal with respect to $g$ if $|f|^n<|g|$ for all $n\in \mathbb{N}$. A positive element $f$ is called {\em dominant} if it is not infinitesimal with respect to any other element. If $g$ is positive and $f$ is infinitesimal with respect to $g$ then we will write $f<<g$. If neither of $f$ and $g$ is infinitesimal with respect to the other one, then we will say that $f$ and $g$ are {\em comparable.} A left-ordered group is called Archimedean if any two non-identity elements are comparable. By a theorem of H\"older, Archimedean groups are Abelian (and hence bi-orderable as well).
 
 \medskip
 
 For a positive element $g\in \Gamma $, we let $\Gamma _g = \{x\in \Gamma : x <<g\}$. For an arbitrary left-orderable group $\Gamma $, the set $\Gamma _g$ may not be very rich with structures, however for bi-orderable groups we have the following easy proposition. 
 
 \begin{prop} \label{thm:subgroup} If $\Gamma $ is bi-orderable and $g$ is a positive element then $\Gamma _g$ is subgroup of $\Gamma $.       \end{prop}
 
 {\bf Proof.} Let $x, y\in \Gamma _g$. Then $|xy|\leq \max\{|x|^2, |y|^2\}$. Consequently, $|xy|^n\leq \max \{|x|^{2n}, |y|^{2n}\}$. Then $xy\in \Gamma _g$. Thus the product of any two elements of $\Gamma _g$ still belongs to it. On the other hand, for an arbitrary $x\in \Gamma _g$, we have $|x^{-1}| = |x|$ hence $x^{-1}\in \Gamma _g$. $\square $
 
 \medskip 
 
 It is easy to see that for an arbitrary element $g$, the subgroup $\Gamma _g$ is not necessarily normal. However, for finitely generated groups we have the following very useful result.
 
 \begin{prop} \label{thm:normalsubgroup} Let $\Gamma $ be a finitely generated bi-orderable group, $S$ be a finite symmetric generating set of $\Gamma $, and $g\in S$ be the biggest element of $S$. Then $\Gamma _g$ is a normal subgroup. Moreover, if $S'$ is another finite symmetric generating set of $\Gamma $ and $g'$ is the biggest element of $S'$ then $\Gamma _g = \Gamma _{g'}$.         
 \end{prop}
 
 {\bf Proof.} Let $x\in \Gamma _g$ be a positive element, and $u\in \Gamma $. Then $|u| < g^m$ for some $m\geq 1$. 
 
 Assuming the opposite, let $uxu^{-1}\notin \Gamma _g$. Then there exists a positive integer $n$ such that $ux^nu^{-1} > g^{2m+1}$. Then $x^n > g$. But then $x\notin \Gamma _g$. Contradiction. 
 
 To prove the second part of the claim, it suffices to observe that there exists $N\geq 1$ such that $|g'| < |g|^N$ and $|g| < |g'|^N$. $\square $ 
 
 \medskip
 
 By the same proof, Proposition \ref{thm:normalsubgroup} can be generalized to infinitely generated groups where the generating set contains a dominant element. 
 
 \begin{prop} \label{thm:normalsubgroup2} Let $\Gamma $ be a bi-orderable group, $S$ be a symmetric generating set of $\Gamma $, and $g\in S$ be a dominant element. Then $\Gamma _g$ is a normal subgroup. Moreover, if $S'$ is another symmetric generating set with a dominant $g'$ then $\Gamma _g = \Gamma _{g'}$.         
 \end{prop}
 
 Proposition \ref{thm:normalsubgroup}  (Proposition \ref{thm:normalsubgroup2}) shows that the subgroup $\Gamma _g$ does not depend on the choice of the finite generating set (of the generating set which contains a dominant element). We will call this subgroup {\em the subgroup of infinitesimals} of $\Gamma $, and denote it by $\mathcal{I}(\Gamma )$. Let us also observe that Proposition \ref{thm:normalsubgroup} is a special case of  Proposition \ref{thm:normalsubgroup2} since any finite symmetric generating set contains a dominant element. Moreover, if a bi-ordered group contains a dominant element then any symmetric generating set also contains a dominant element. The class of bi-orderable groups with a dominant element contains all finitely generated bi-orderable groups.
 
 \medskip
  
 We would like to make more observations on the properties of the subgroup of infinitesimals. All these propositions will be stated for bi-ordered groups with a dominant element. The following proposition follows from H\"older's Theorem which states that an Archimedean group is Abelian. 
 
 \medskip
 
 \begin{prop} \label{thm:abelian} Let $\Gamma $ be a bi-ordered group with a dominant element. Then the quotient $\Gamma /\mathcal{I}(\Gamma )$ is Abelian.
 \end{prop}
 
 \medskip
 
 \begin{cor} \label{thm:contains} As a corollary we also obtain that $\mathcal{I}(\Gamma )\supseteq [\Gamma , \Gamma ]$. 
 \end{cor} 
 
 \medskip
 
 The above observations (i.e. Propositions \ref{thm:subgroup}, \ref{thm:normalsubgroup}, \ref{thm:abelian} and Corollary \ref{thm:contains}) have been made in \cite{A} but only for a finitely generated groups and also in a very special setting. Namely, let $\Gamma \leq \mathrm{Homeo}_{+}(I)$ such that every non-identity element has finitely many fixed points. Then we can introduce a bi-order in $\Gamma $ as follows: for all $f, g\in \Gamma $ we say $f < g$ if $f(x) < g(x)$ near zero, i.e. in some interval $(0, \epsilon )$. For this particular bi-order, the first named author has observed in \cite{A} that the claims of the above propositions hold, in fact, these observations have been much exploited in \cite{A}. We also would like to remark that the concept of infinitesimals has been used by Farb-Franks in \cite{FF} and this has inspired the first named author to study its properties further.
 
 \medskip 
 
  Let us also observe the following nice property of the subgroup of infinitesimals.
  
  \begin{prop} \label{thm:invariance} Let $\Gamma $ be a bi-ordered group with a dominant element, and $\phi : \Gamma \to \Gamma $ be an automorphism of $\Gamma $ which preserves a bi-order in $\Gamma $. Then the restriction $\phi |_{\mathcal{I}(\Gamma )}$ is an automorphism of $\mathcal{I}(\Gamma )$. $\square $
  \end{prop}  
  
  \medskip 
  
 Let us recall that an automorphism is said to preserve a left order if positive elements are mapped to positive elements. Such a property for an automorphism turns out to be quite useful in many situations. For example, if a bi-ordered group $H$ acts on a bi-ordered group $G$ with order preserving  automorphisms then one can put a bi-order on the semi-direct product $H\ltimes G$ as follows: let $\prec _1, \prec _2$ be the bi-orders in the groups $H, G$ respectively. An element of $H\ltimes G$ can be written as a pair $(h,g)$. We let $(h,g)$ be positive iff either $1\prec _1 h$ or $h = 1$ and $1\prec _2 g$. Since the action of $H$ preserves the bi-order $\prec _2$ on $G$, one can check easily that the bi-order on $H\ltimes G$ is well defined.
  
  \medskip
  
 Now we would like to prove a modification of the above propositions which will have an important application to knot groups. We will state this proposition for semidirect products of $\mathbb{Z}\ltimes G$, although it can be stated for proper HNN extensions as well (the knot groups that we have in mind do indeed split as semidirect products).
 
 \medskip
 
 \begin{prop} \label{thm:semidirect} Let $G$ be a bi-ordered group with a dominant element, $\phi :G\to G$ be an automorphism of $G$ which preserves the bi-order of $G$, and $\mathbb{Z}\ltimes G$ be the semidirect product corresponding to the action of $\phi $. Then the groups $\mathbb{Z}\ltimes \mathcal{I}(G)$ and $\mathbb{Z}\ltimes (G/\mathcal{I}(G))$ are bi-orderable. In particular, if any two elements of $G\backslash [G,G]$ are comparable then $\mathbb{Z}\ltimes (G/[G,G])$ is bi-orderable.
 \end{prop}  
 
 {\bf Proof.} Let $<$ be a bi-order in $G$. The restriction of $\phi $ to $\mathcal{I}(G)$ preserves the bi-order $<$ in it. Then the bi-orderability of the group $\mathbb{Z}\ltimes \mathcal{I}(G)$ follows from Proposition \ref{thm:invariance}. 
 
 \medskip
 
 For the bi-orderability of the second group $\mathbb{Z}\ltimes (G/\mathcal{I}(G))$, we define a bi-order $\prec $ in $G/\mathcal{I}(G)$ as follows. For any two elements $\xi , \zeta \in G/\mathcal{I}(G)$ let $x, y$ be the representatives of $\xi , \zeta $ respectively. We let $\xi \prec \zeta $ if and only if $x < y$. It is an easy check that $\prec $ is indeed a well defined bi-order in $G/\mathcal{I}(G)$. 
 
 \medskip
 
 If any two elements of $G/[G,G]$ are comparable, by Corollary \ref{thm:contains} we obtain that $\mathcal{I}(G) = [G,G]$. Therefore $\prec $ is a well defined bi-order in $G/[G,G]$. Moreover, $\prec $ is preserved by $\phi $. $\square $ 
 
 \medskip
 
 \begin{prop} \label{thm:eigenvalues} Let $\mathbb{Z}$ act on $\mathbb{Z}^d$ by an automorphism $A\in GL(d, \mathbb{Q})$.
 
 a) If $\mathbb{Z}\ltimes _{A}\mathbb{Z}^d$ is bi-orderable then $A$ has at least one positive eigenvalue. 
 
 b) Let $\Lambda $ be the set of all eigenvalues of $A$ which are either imaginary or negative, and $V_{\Lambda }$ be the eigenspace corresponding to $\Lambda $. If $\mathbb{Z}\ltimes _{A}\mathbb{Z}^d$ is bi-orderable then $V_{\Lambda }$ does not contain a non-zero rational vector.

 \end{prop}  
 
 {\bf Proof.} In Proposition 4.16 in \cite{CR} it is proved that the if all eigenvalues of $A$ are positive and real then the semidirect product $\mathbb{Z}\ltimes \mathbb{R}^d$ is bi-orderable. Part a) is a partial converse to this claim, and it might also follow from the arguments of \cite{PR}, \cite{CR}. We will nevertheless prove it for the case of actions on $\mathbb{Z}^d$ (instead of $\mathbb{R}^d$). Somewhat similar ideas have been used also in the proof of Proposition 2.2. in \cite{LRR}.
 
 \medskip
 
  Let $\prec $ be a bi-order in $\mathbb{Z}\ltimes _{A}\mathbb{Z}^d$, and $x_1, \dots , x_d$ be positive generators of $\mathbb{Z}^d$. By the bi-orderability, the action of $A$ on $\mathbb{Z}^d$ preserves the bi-order $\prec $ on $\mathbb{Z}^d$. Without loss of generality, we may also assume that all elements of $\mathbb{Z}^d\backslash \{0\}$ are comparable (otherwise we will concentrate on the semidirect product $\mathbb{Z}\ltimes _{A}\mathcal{I}(\mathbb{Z}^d)$ and proceed by induction; notice that $\mathcal{I}(\mathbb{Z}^d)\cong \mathbb{Z}^r$ where $r < d$). Then there exists a unique vector $u = (\alpha _1, \dots , \alpha _d)$ up to a scalar multiplication such that all entries are positive and for all $i, j\in \{1, \dots , d\}$, $x_i^{m} \prec x_j^n$ iff $\frac{m}{n} < \frac{\alpha _j}{\alpha _i}$.  Thus we obtain that for all $z = (z_1, \dots , z_d)\in \mathbb{Z}^d$ we have $0 \prec z$ only if $z_1\alpha _1 + \dots + z_d\alpha _d \geq 0$. In other words, if $0 \prec z$ then $\langle z, u \rangle \geq 0$ where $\langle ., . \rangle $ denotes the standard dot product in $\mathbb{R}^d$.   
   
   \medskip
   
  Now, assume that $A$ has a non-positive real eigenvalue $\lambda $. Then $\lambda $ is negative. Let $v = (v_1, \dots , v_d)\in \mathbb{R}^d$ be an eigenvector of $A$ with respect to $\lambda $. If the vector $v$ is projectively rational (i.e. it is a scalar multiple of a vector with all coordinates rational) then $\lambda $ is rational, and we obtain that $A(Mw) = -Nw$ for some integral vector $w$ and natural numbers $M$ and $N$. This violates bi-orderability of $\mathbb{Z}\ltimes _{A}\mathbb{Z}^d$ since the action of $A$ preserves $\prec $; thus we can assume that the vector $v$ is not projectively rational. Then for every $\epsilon > 0$ there exists a vector $z = (z_1, \dots , z_d)\in \mathbb{Z}^d$ and a scalar $c > \frac{1}{\epsilon }$ such that $0 < ||cv-z|| < \frac{\epsilon }{||A||+1}$ {\bf(1)}. Then $0 < ||cAv-Az|| < \epsilon $ {\bf(2)}.  Without loss of generality, we may assume that $0 \prec z$. Then $\langle z, u \rangle \geq 0$, i.e. the angle between $z$ and $u$ is not obtuse. If $\langle v, u \rangle > 0$ (or if $\langle v, u \rangle < 0$) then for sufficiently small $\epsilon $, the inequalities {\bf (1)} and {\bf (2)} imply that $z$ and $Az$ cannot be both positive which is again a contradiction. Thus $\langle v, u \rangle = 0$, i.e. an eigenvector corresponding to a negative eigenvalue is orthogonal to $u$.    
   
   \medskip
   
   Similarly, in the case of an imaginary eigenvalue, we obtain that the corresponding eigenvector must be orthogonal to $u$ (we leave this as an exercise to the reader). 
   
   \medskip
   
   Thus the eigenspace $V_{\Lambda }$ is orthogonal to $u$. Now part a) follows immediately (since $u$ is a non-zero vector, $V_{\Lambda }$ is a proper subspace of $\mathbb{R}^d$). If  $V_{\Lambda }$ contains a non-zero rational vector $w$, let $V(w)$ be the minimal $A$-invariant $\mathbb{Q}$-subspace of $V$ containing $w$, and $V_{\mathbb{Z}}(w)$ be the subgroup of all integral vectors of $V(w)$. Then $\mathbb{Z}\ltimes V_{\mathbb{Z}}(w)$ is bi-orderable. But this contradicts part a) since the acting automorphism, when restricted to $V_{\Lambda }$, has no positive eigenvalue.   $\square $
   
   \medskip
   
   \begin{rem}\label{thm:biorderabilityremark} Indeed, we also prove that if $A\in GL(d,\mathbb{Q})$ is an integral matrix then a semidirect product $\mathbb{Z}\ltimes _{A}\mathbb{Z}^d$ is bi-orderable provided $A$ has at least one positive eigenvalue, and no $(d-1)$-dimensional eigenspace of $A$ contains a non-zero rational vector. 
   \end{rem}
   
   To point out further immediate corollaries of the proof, let us recall the  Rational Form Theorem that for any matrix $A\in GL(d, \mathbb{Q})$, there exists a decomposition  $\mathbb{Q}^d = V_1\times \dots \times V_r$ into irreducible components corresponding to the factorization of the characteristic polynomial of $A$. Thus the matrix $A$ has a decomposition into $r$ blocks corresponding to the restriction of the map $A$ to the subspaces $V_1, \dots , V_r$. A vector $x\in V$ can be written uniquely as $x = x_1 + \dots + x_k$ where $k\leq r$, $x_i \in V_{n_i}$ for all $1\leq i\leq k$ and $n_1, \dots , n_k$ are all distinct. We will call the vector $x_i$ and irreducible component of $x$. The matrix $A$ can be canonically written as a matrix of $r$ irreducible blocks corresponding to the subspaces $V_1, \dots , V_r$.
   
   \medskip  
   
    Now we observe the following immediate but very important corollary of Proposition \ref{thm:eigenvalues} and Rational Form Theorem which provides a classification of all bi-orderable semidirect products $\mathbb{Z}\ltimes \mathbb{Z}^n$.
   
   \begin{thm}\label{thm:classification} Let $\mathbb{Z}$ act on $\mathbb{Z}^d$ by an automorphism $A\in GL(d, \mathbb{Q})$. Then the group $\mathbb{Z}\ltimes _{A}\mathbb{Z}^d$ is bi-orderable if and only if all irreducible blocks of the matrix $A$ over the rationals possesses a positive real eigenvalue.
   \end{thm} 
   
   Proposition \ref{thm:eigenvalues} and Rational Form Theorem have another immediate corollary which we state in the following very useful lemma
   
   \begin{lem}\label{thm:use} Let $V$ be a bi-ordered Abelian group, $T:V\to V$ be an order-preserving automorphism, $V_0\leq V$ be a $T$-invariant subgroup of finite index such that $V/V_0 \cong \Z ^d$ for some $d\geq 1$, and $x\in V\backslash V_0$ such that $x$ is comparable to $Tx$. Let also $A\in GL(d,\mathbb{Q})$ be the matrix representing the map $\overline{T}:V/V_0\to V/V_0$ induced by $T$, and $\overline{x}\in V/V_0$ be the image of $x$ in the quotient $V/V_0$. If $A$ consists of a single irreducible block over the rationals then the restriction of $\overline{T}:V/V_0\to V/V_0$ possesses a positive real eigenvalue.
   \end{lem}

 \bigskip
 
 In the remainder of this section, we would like to introduce the following very useful notion of a completion of a free bi-ordered group.  Let $G$ be a free group with a bi-order $<$. We'll introduce a completion $\overline{G}$ of $G$. Our main goal is to extend the group $G$ to include ``the limits" of sequences of type $f^ngf^{-n}, n\geq 1$ (notice that these sequences are always monotone and bounded). Inductively, we will build a sequence of linearly ordered sets $G_n, n \geq 0$ such that

\medskip
 
 i) $G_0 = G$;

\medskip
 
 ii) $G_{n}\subseteq G_{n+1}$, for all $n\geq 0$;

\medskip
 
 iii) the linear order on $G_0$ coincides with $<$; 

\medskip
 
 iv) for all $n\geq 0$, the linear order on $G_{n+1}$ is an extension of the linear order on $G_n$;
 
  \medskip
  
 v) for all $n\geq 0$, a multiplication of some elements of $G_n$ is defined in $G_n$, moreover, in $G_{n+3}$, a multiplication between any two elements of $G_n$ is defined, i.e. $G_nG_n\subseteq G_{n+3}$.
  
 \medskip
 
 vi) for all $n\geq 0$, a multiplication in $G_{n+1}$ is an extension of the multiplication in $G_n$. 
 
 \medskip
 
 Indeed, we will use the following  very specific extension of the group operation and the linear order: let $G_0 = G$, $S^{(0)}$ be the set of all sequences in $G_0$ of the form  $(f^ngf^{-n})_{n\geq 1}$. Notice that all sequences in $S^{(0)}$ are monotone and bounded, moreover, a constant sequence $(g, g, \dots )$ belongs to $S^{(0)}$ for all $g\in G_0$. Thus we may view elements of $G_0$ as constant sequences hence as elements of $S^{(0)}$. 
 
 \medskip
 
 We will build a sequence of sets $A_n, n\geq 0$ such that for all $n\geq 1$ a multiplication between any two element of $A_{n-1}$ is defined and the following conditions hold: 
 
 1) for all $n\geq 0$, $A_n = \Omega _n\sqcup \Delta _n$ where $\Omega _n$ is a set of infinite sequences of elements from $A_{n-1}$ for all $n\geq 1$, and $\Delta _n$ is a set of infinite sequences of elements from $\Omega _n$ for all $n\geq 0$.
 
 2) for all $n\geq 1$, $A_{n-1}$ is a subset of $\Omega _n$ where every element of $A_{n-1}$ is viewed as a constant sequence in $\Omega _n$. 
 
  3) for all $n\geq 0$, $\Omega _{n}$ is a subset of $\Delta _n$ where every element of $\Omega _n$ is viewed as a constant sequence in $\Delta _n$. 
 
 4) for all $n\geq 1$, $A_{n-1}A_{n-1} \subseteq A_n$.
 
 \medskip
  
 We let $\Omega _ 0 = G_0, \Delta _0 = S^{(0)},$ and $A_0 = \Omega _0 \sqcup \Delta _0$.  Assume now the sets $A_i, 0\leq i\leq n$  are all defined such that conditions 1)-4) hold up to the index $n$.  We will define $A_{n+1}$ as follows:
 
 \medskip
 
 Let $$\Omega ^{+} = \{(fg_k)_{k\geq 1} \ | \ f\in A_{n-1}, (g_k)_{k\geq 1} \ \mathrm{belongs \ to} \ A_n\}, \ \mathrm{and}$$ \   $$\Omega ^{-} = \{(g_kf)_{k\geq 1} \ | \ f\in A_{n-1}, (g_k)_{k\geq 1} \  \mathrm{belongs \ to} \ A_n\}.$$ Notice that since $f\in A_{n-1}$ and $g_k\in A_{n-1}, \forall k\geq 1$, by condition 4), the multiplications are already defined in $A_n$. We let $\Omega _{n+1} = \Omega ^{+}\cup \Omega ^{-}$.  
 
 \medskip
 
  For any two sequences $a = (a_k)_{k\geq 1}$ and $b = (b_k)_{k\geq 1}$ in $A_{n}$, for all $m\geq 1$,  we define $ab_m$ as a sequence $(a_kb_m)_{k\geq 1}$ (thus $ab_m$ belongs to $\Omega _{n+1}$). Then we define $ab$ as the sequence $(ab_k)_{k\geq 1}$, and $a\star b$ as the sequence $(a_kb)_{k\geq 1}$. Then  we let $$\Delta _{n+1} = \{ab \ | \ a, b\in A_n\}\cup \{a\star b \ | \ a, b\in A_n\} \ \mathrm{and} \  A_{n+1} = \Omega _{n+1}\sqcup \Delta _{n+1}.$$ Notice that on the set $A_{n+1}$ all conditions 1)-4) still hold thus we may proceed by induction.

 \medskip
 
 We now let $B_{2n} = \Omega _n, B_{2n+1} = \Delta _n$ for all $n\geq 0$. On the sets $B_n, n\geq 0$ we will introduce equivalence relations and denote the set of equivalence classes of $B_n$ by $G_n$. Then on $G_n$ we will introduce a linear order. We do this as follows:
 
 \medskip
 
 First, on $G_0 = B_0 = G, G_1 = B_1$ where the bi-order on $G_0$ coincides with the bi-order $<$ on $B_0 = G$, and the bi-order on $G_1$ is defined as follows: $G_1 = B_1 = \Delta _0$ consists of bounded monotone sequences of the form $(f^ngf^{-n})_{n\geq 1}$; we let this sequence be positive iff $g$ is positive.
 
 \medskip
 
 Assume now the equivalence classes $G_0, \dots , G_{n-1}$ of $B_0, \dots , B_{n-1}$ are all defined for some $n\geq 2$. Then on the set $B_i, 2\leq i\leq n-1$ we can also define a partial order (by abuse of notation we still denote it by $<$) by letting $x<y$ for $x, y\in B_i$ iff $[x]<[y]$ where $[x], [y]\in G_i$ represent equivalences classes of elements in $x$ and $y$ respectively.
 
 \medskip
 
 Then on the set $B_{n}, n\geq 1$ we introduce the following equivalence relation: the elements $\bar{f} = (f_n)_{n\geq 1}, \bar{g} = (g_n)_{n\geq 1}$ of $B_n$ (which are sequences in $B_{n-1}$) are called equivalent if the following conditions hold:
 
 \medskip
 
 1) if either both $\bar{f} = (f_n)_{n\geq 1}, \bar{g} = (g_n)_{n\geq 1}$ are non-decreasing or both are non-increasing then there exists $N\geq 1$ such that for all $k\geq N$ there exists $p, q, r, s$ such that $p < q, r < s$, and $f_p \leq g_k \leq f_q, g_r \leq f_k \leq g_s$. 
 
 \medskip
 
 2) if $\bar{f}$ is non-decreasing and $\bar{g}$ is non-increasing then there exists at most one $\omega \in G_n$ such that $f_i\leq \omega \leq g_j$, for all $i, j\geq 1$. 
 
 \medskip
 
 3) if $\bar{f}$ is non-increasing and $\bar{g}$ is non-decreasing then there exists at most one $\omega \in G_n$ such that $g_i\leq \omega \leq f_j$, for all $i, j\geq 1$.   
 
 \medskip
 
 We have already defined $G_{n}$ as the set of equivalence classes of $B_{n}$. One can check that for all $a, b\in A_n$, the elements $ab$ and $a\star b$ of $A_{n+1}$ will belong to the same equivalence classes thus they represent the same elements of $G_n$. Therefore, for  $a, b\in G_n$, we can denote the multiplication simply by $ab$. In the sequel, by abuse of notation we will still write the elements of $G_i, i\geq 1$ as sequences of elements of $B_{i-1}$ representing these elements; our discussion will not depend on the choices of representatives of equivalence classes.   
 
 \medskip
 
 Now it remains to show the extension of the bi-order $<$ in $G_0\cup G_1$ to all $G_n, n\geq 2$. We will extend it inductively as follows: assume now that the it is defined on all $G_i, 0\leq i\leq n$ for some $n\geq 1$ such that for all $0\leq i\leq n-1$, on $G_{i+1}$ it is an extension of the bi-order on $G_i$, and all elements of $G_{i+1}$ consists of bounded monotone sequences of $G_i$. Inductively, we extend the order $<$ on $G_{n}$ to $G_{n+1}$ as follows: let $\bar{f} = (f_n)_{n\geq 1}, \bar{g} = (g_n)_{n\geq 1}$ be elements of $G_n$ as monotone sequences in $G_{n-1}$. We say  $\bar{f} < \bar{g}$ in one of the following cases:
 
 i) both $\bar{f}, \bar{g}$ are non-decreasing and there exists $k\geq 1$ such that $g_k > f_n$ for all $n\geq 1$;
 
 ii) both $\bar{f}, \bar{g}$ are non-increasing and there exists $k\geq 1$ such that $f_k < g_n$ for all $n\geq 1$;
 
 iii) $\bar{f}$ is non-decreasing, $\bar{g}$ is non-increasing, and  $f_i < g_j$ for all $i, j\geq 1$;
 
 iv) $\bar{f}$ is non-increasing, $\bar{g}$ is non-decreasing, and there exists $k\geq 1$ such that $f_k < g_k$.
 
 \medskip
 
 By the way we defined the sets $A_n, n\geq 0$, all elements of $G_{n+1}$ will consists of bounded monotone sequences thus we may proceed by induction. 
 
 \medskip
 
  Thus both the multiplication and the the order $<$ is extended to the whole union $\cup _{n\geq 0}G_n$. We aslo define the inverse of an element $\bar{f} = (f_n)_{n\geq 1}$ inductively as $\bar{f}^{-1} = (f_n^{-1})_{n\geq 1}$. One can check that $\cup _{n\geq 0}G_n$ becomes a bi-ordered group. We denote this group by $G^{(1)}$ (let us also denote $G$ by $G^{(0)}$). Then we take a set $S^{(1)}$ of all sequences in $G^{(1)}$ of the form  $(f^ngf^{-n})_{n\geq 1}$, and repeat the process again to obtain a bi-ordered group $G^{(2)}$, and so on. The union $\cup _{n\geq 0}G^{(n)}$ becomes a bi-ordered group. We denote it with $\overline{G}$, and call it {\em the completion of} $G$. 
 
 \medskip
 
 We would like to end this section with the following useful notion.
 
 \begin{defn} Let $G$ be a group with a bi-order $<$. Two elements $f, g\in G$ are called weakly comparable if there exists $h\in G$ such that $f$ and $ghg^{-1}$ are comparable.
 \end{defn}
 
 Often times an automorphism of a bi-orderable group sends an element of a group to another element which is not comparable to it. This turns out to be a serious issue in the situation interesting to us; fortunately, we will be able to establish a weak comparability which will be sufficient for our purposes.

 \bigskip
 
 \section{Non-bi-orderability of $\pi _1(6_2)$}
 
 \medskip
 
 Let $K\subset \mathbb{S}^3$ represent the knot $6_2$. In Figure 1, we have shown the picture of the Seifert surface $\Sigma $ of the knot $6_2$. Let $K$ denote the boundary of this surface. In Figure 3, we have chosen the base point for $\pi _1(\mathbb{S}^3\backslash K)$ and make some convenient choices for the generators of $\pi _1(\mathbb{S}^3\backslash K)$ and $\pi _1(\Sigma \times \{1\})$. In the figure, the loops $\alpha , \beta , \gamma , \delta $ generate $\pi _1(\Sigma \times \{1\})$, and the loops $x, a, b, c$ generate $\pi _1(\mathbb{S}^3\backslash K)$.
 
 \medskip
 
 We obtain a presentation $$\langle t, \alpha , \beta , \gamma , \delta \ | \ t\alpha t^{-1} = \alpha _{-}, t\beta t^{-1} = \beta _{-}, t\gamma t^{-1} = \gamma _{-}, t\delta t^{-1} = \delta _{-}\rangle $$ where $\alpha _{-}, \beta _{-}, \gamma _{-}, \delta _{-}$ are the loops on $\Sigma \times \{-1\}$ parallel to $\alpha , \beta , \gamma , \delta $ respectively. By computing these loops we obtain the following presentation for the knot group $6_2$: $$\Gamma = \langle t, x, a, b, c \ | \ ta^{-1}t^{-1} = xb, txat^{-1} = x, tbt^{-1} = c^{-1}, tct^{-1} = abc\rangle $$
 
 \medskip
 
 Both of the subgroups generated by $a^{-1}, xa, b, c$ and $xb, x, c^{-1}, abc$ coincide with the free group generated by $x, a, b, c$ hence the group $\Gamma $ by the above presentation can be given as a semidirect product $\mathbb{Z}\ltimes \mathbb{F}^4$ where the acting group $\mathbb{Z}$ is generated by $t$, and the free group $\mathbb{F}_4$ is generated by $x, a, b, c$. 
 
 \medskip
 
 \begin{figure}[h!]
  \includegraphics[width=6in,height=5in]{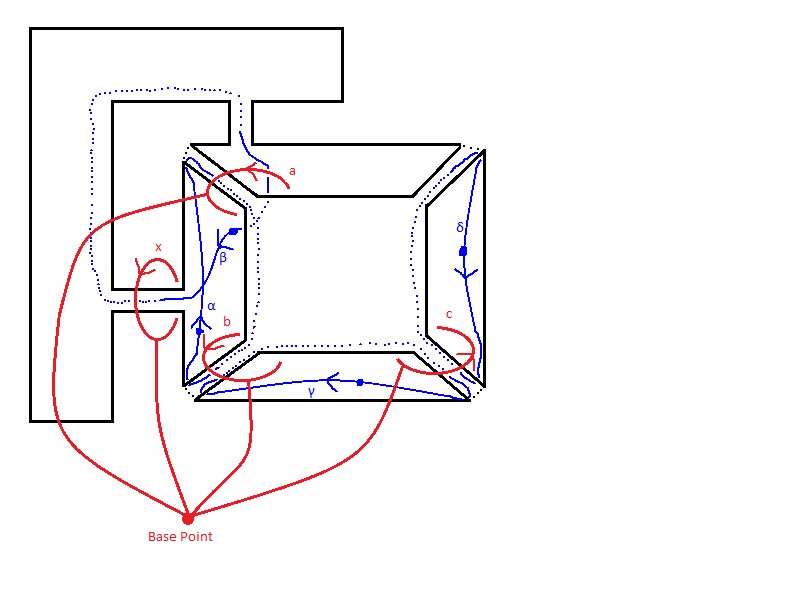}
\caption{A Seifert surface of $6_2$. Blue loops represent the generators of the fundamental group of the surface, while the red loops represent the generators of the fundamental group of complement of the surface.}
\end{figure}
 
 \medskip
 
 Let $G$ be a subgroup of $\Gamma $ generated by $x, a, b$ and $c$ (thus $G$ is isomorphic to $\mathbb{F}_4$). Let us assume that $\Gamma $ has a bi-order $<$. Then this bi-order can be restricted to $G$ and the conjugation by $t$ preserves this bi-order in $G$. Notice that the quotient $\Gamma /[G,G]$ is isomorphic to $\mathbb{Z}\ltimes _{M}\mathbb{Z}^4$ where the automorphism is given by the matrix \[  M =  \left ( \begin{tabular}{cccc} 2 &  1 & 0 & 0 \\  0 & 0 & 0 & 1 \\  1 & 1 & 0 & -1 \\ 0 & 0 & 1 & 1  \end{tabular}   \right ) \]
 
 \medskip
 
 The characteristic polynomial det$(M - \lambda I)$ of this matrix equals $f(\lambda ) = \lambda ^4 - 3\lambda ^3 + 3\lambda ^2 - 3\lambda + 1$ which has exactly two positive real roots and no rational roots.   \footnote{It is not accidental that this polynomial coincides with the Alexander polynomial of the knot. One can involve the rational root test to show that it does not have a rational root. To see why it does not have four real roots, first, notice that the second derivative $f''(\lambda ) = 12\lambda ^2 - 18\lambda + 6$ has exactly two roots: $\lambda _1 = \frac{1}{2}$ and $\lambda _2 = 1$, and it is negative in the interval $(\frac{1}{2}, 1)$ . On the other hand, in this interval, $f(\lambda ) = \lambda ^4 - \lambda ^3 + 2(\lambda ^2 - \lambda ^3) + (\lambda ^2 - \lambda ) -2\lambda + 1 < 0$. If $f(\lambda )$ has four real roots then, since $f(0) > 0$, it has either three roots in $(0, \frac{1}{2})$ and one root in $(1,\infty )$, or it has one root in $(0, \frac{1}{2})$ and three roots in $(1,\infty )$. By Rolle's Theorem, one can deduce that $f''$ must have a root in the interval $(0, \frac{1}{2})\cup (1,\infty )$. Contradiction.} 
 
 \medskip
 
 Since $f(\lambda )$ has no rational root, any two elements of $G\backslash [G, G]$ are comparable. Then by Proposition \ref{thm:semidirect} the quotient $\Gamma /[G,G]$ must be bi-orderable. The conditions of Proposition \ref{thm:eigenvalues} are not met, and we do not obtain a contradiction at this level; since $M$ does not have a rational eigenvalue, by Remark \ref{thm:biorderabilityremark}, the quotient $\Gamma /[G, G]$ is indeed bi-orderable. 
 
 \medskip
 
 Now we consider the central series $G_n, n = 0, 1, 2, \dots $ of $G$ defined as follows: $G_0 = G, G_{n+1} = [G_n,G_n]$, for all $n\geq 0$. Each of the subgroups $G_n$ is invariant under the conjugation by $t$. Thus $t$ acts on all quotients $G_n/G_{n+1}$. We have already studied this action for $n=0$; now, let us study it for $n=1$.  Let $$u = [x,a], v = [x,b], w = [x,c], \xi = [a,b], y = [a,c], z = [b,c].$$

 \medskip
 
 By abusing the notation we will denote the images of the elements $u, v, w, \xi , y, z$ in the quotient $G_1/G_2$ by the same letters. In the group $G_1/G_2$ we have a conjugation by $t$ given as $$tut^{-1} = v^{-1}, \ tvt^{-1} = w^{-2}z^{-1}, \ twt^{-1} = u^2v^2w^2\xi ^{-1}z,$$ \ $$t\xi t^{-1} = wz, \ tyt^{-1} = u^{-1}v^{-1}w^{-1}\xi z^{-1}, \ tzt^{-1} = yz,$$ thus the $t$-conjugation acts by an automorphism given by a matrix \begin{displaymath}  N =  \left ( \begin{array}{cccccc} 0 &  -1 & 0 & 0 & 0 & 0 \\  0 & 0 & -2 & 0 & 0 & -1 \\  2 & 2 & 2 & -1 & 0 & 1 \\ 0 & 0 & 1 & 0 & 0 & 1 \\ -1 & -1 & -1 & 1 & 0 & -1 \\ 0 & 0 & 0 & 0 & 1 & 1  \end{array}   \right ). \end{displaymath}
 
 \medskip

 The characteristic polynomial of the matrix $N$ equals $p(\lambda ) = \lambda ^6 - 3\lambda ^5 +8\lambda ^4 -12\lambda ^3 + 8\lambda ^2 - 3\lambda +1$. It has a double root $\lambda _1 = \lambda _2 = 1$, and 4 other imaginary roots $\lambda _3, \lambda _4, \lambda _5$ and $\lambda _6$. Let $E_1$ be the two-dimensional real eigenspace corresponding to the eigenvalue 1, and $E_2$ be the 4-dimensional real  eigenspace corresponding to the eigenvalues $\lambda _3, \lambda _4, \lambda _5$ and $\lambda _6$. The subspace $E_2$ contains a non-zero rational vector, moreover, it does not contain a proper $A$-invariant subspace $V$ such that $V$ contains a non-zero rational vector. Thus if the semidirect product $\mathbb{Z}\ltimes (G_1/G_2)$ was bi-orderable where the action is given by the $t$-conjugation then we would obtain a contradiction by Theorem \ref{thm:classification}. However, the bi-order $<$ does not necessarily induce a bi-order on $G_1/G_2$ which is preserved by the $t$-conjugation. 
 
 \medskip
 
 To resolve this issue we would like to consider two cases: 
 
 \medskip
 
 {\em Generic Case:} Let us first assume that none of the elements $u, v, w, \xi , y, z$ are weakly comparable.

 Let us observe that in the group $G = G_0$, for any $\theta \in \{u, v, w, \xi , y, z\}$ and $n\geq 1$, we have $t^n\theta t^{-n} = g_{1,n}\omega _1 g_{1,n}^{-1}\dots g_{k,n}\omega _k g_{k,n}^{-1}$ where $$\omega _i\in \{u, u^{-1}, v, v^{-1}, w, w^{-1}, \xi , \xi^{-1}, y, y^{-1}, z, z^{-1}\}$$ and $g_{i,n}\in G$ for all $1\leq i\leq k$ (here, $k$ depends on $n$, so we will also write $k = k(n)$). The conjugate $g_{j,n}\omega _j g_{j,n}^{-1}, 1\leq j\leq k(n)$ will be called {\em a term of $t^n\theta t^{-n}$ corresponding to $\omega _j$}. We will call an element of $G_1$ {\em admissible} if it can be represented by a reduced non-trivial word $W$ in the alphabet $$\{u, u^{-1}, v, v^{-1}, w, w^{-1}, \xi , \xi^{-1}, y, y^{-1}, z, z^{-1}\}$$ such that for all $\omega _i$ in the alphabet, if $tWt^{-1}$ has a term corresponding to $\omega _i$ then it has no term corresponding to $\omega _i^{-1}$. 
 
 \medskip
 
 Now, let us make an observation that, for all $\omega _i$ from the alphabet, if $t\theta t^{-1}$ has a term corresponding to $\omega _i$ then it has no term corresponding to $\omega _i^{-1}$, i.e. each of letters $\omega _i$ are admissible elements. \footnote{It follows by the action of $t$ on $G_1/G_2$ that, in the group $G = G_0$, the element $tut^{-1}$ will have only one term, and it corresponds to $v^{-1}$; the element $tvt^{-1}$ has two terms corresponding to $w^{-1}$ and one term corresponding to $z^{-1}$; the element $twt^{-1}$ has two terms corresponding to each of $u, v, w$, one term corresponding to $\xi ^{-1}$, and one term corresponding to $z$, and so on.}  Recall that two elements $f,g\in G$ are called weakly comparable if $g$ and $hfh^{-1}$ are comparable for some $h\in G$. 
 
 \medskip
 
 Considering any seven different positive values of the exponent $m$, by pigeonhole principle, we can claim that for some $$\omega \in \{u, u^{-1}, v, v^{-1}, w, w^{-1}, \xi , \xi^{-1}, y, y^{-1}, z, z^{-1}\},$$ and for two of those seven values of $m$, there exists $i\in \{1, \dots , k(m)\}$ such that  all terms of $t^m\omega t^{-m}$ not corresponding to $\omega $ is infinitesimal with respect to the term $g_{i,m}\omega g_{i,m}^{-1}$. Let $p, q$ be these two desired values of the exponent $m$. Then, by replacing the $t$-conjugation with a $t^rg $-conjugation if necessary, for some $g \in G$ (more precisely, for $g = g_{i,p}^{-1}g_{i,q}$), we may and will assume that all the conjugates $t^{n}\omega t^{-n}, n\in \mathbb{Z}$ are comparable where $r = |p-q|\leq 6$. 
 
 \medskip
 
 In the quotient $G_1/G_2\cong \Z ^6$, the images of all the conjugates $t^{n}\omega t^{-n}, n\in \mathbb{Z}$ generate the entire quotient $G_1/G_2$. Then we can find an element $\eta = \eta (u, v, w, \xi, y, z)\in G_1$ such that, in the quotient $G_1/G_2$, the image of $\eta $ belongs to the 4-dimensional irreducible component isomorphic to $\Z ^4$.    
 
 \medskip

 Let $H$ be the subgroup of $G_1$ generated by the subset $S = \{t^{n}\eta t^{-n} \ | \ n\in \mathbb{Z}\}$ of all conjugates of $\eta $. All elements of $S$ are comparable to each other. Let us observe that a conjugation by an element $g \in G$ induces an identity map in the quotient $G_1/G_2$. Then all elements of $[H,H]$ are infinitesimal with respect to all elements of $S$, and we obtain that the semidirect product $\mathbb{Z}\ltimes (H/\mathcal{I}(H))$ is bi-orderable, where the action is given by the  $t$-conjugation. By Proposition \ref {thm:abelian} the group $H/\mathcal{I}(H)$ is Abelian. Let $L = H\cap G_2$. Notice that $L$ is invariant under the $t$-conjugation, and $H/L\cong \Z ^4$. Let $T:H/\mathcal{I}(H)\to H/\mathcal{I}(H)$ denotes the map induced by the $t$-conjugation. Let also $L'$ be the image of $L$ in $H/\mathcal{I}(H)$.  Then we apply Lemma \ref{thm:use} by taking $V = H/\mathcal{I}(H), V_0 = L'$ and $x = \eta $ to obtain a contradiction.

 \medskip
  
 {\em Non-generic Case:} Some two of the elements $u, v, w, \xi , y, z$ are weakly comparable.
 
 \medskip
  
 To treat this case first we choose a positive element $\eta = \eta (u, v, w, \xi, y, z)\in G_1$ such that, in the quotient $G_1/G_2$, the image of $\eta $ belongs to the 4-dimensional irreducible component isomorphic to $\Z ^4$. Without loss of generality let us assume that $T(\eta ) \leq \eta $. The sequence $(\eta ^n)_{n\geq 1}$ is increasing in $G$. If it is not bounded then $\eta \leq T^{-1}(\eta ) < \eta ^n$ for some $n\geq 1$. Thus $T^{-1}(\eta )$ (hence also $T(\eta )$) is comparable with $\eta $. Then we again let $H$ be the subgroup of $G_1$ generated by the subset $S = \{t^{n}\eta t^{-n} \ | \ n\in \mathbb{Z}\}$ of all conjugates of $\eta $, and proceed as in the Generic Case. 
 
 \medskip
 
 Thus we may assume that $(\eta ^n)_{n\geq 1}$ is bounded. Let us denote this sequence by $\omega $. Although it is not necessarily an element of $\overline{G}$ we can define ``a conjugation" by it which we denote by $I_{\omega }$, namely, for all $x\in \overline{G}$ let $I_{\omega }(x)$ denotes the sequence $(\eta ^nx\eta ^{-n})_{n\geq 1}$. Notice that this sequence is bounded and monotone hence it can be viewed as an element of $\overline{G}$. In the completion $\overline{G}$ of $G$, we will have $\eta (I_{\omega }T(\eta ))\eta ^{-1} =  I_{\omega }T(\eta )$, and more generally,  $\eta (I_{\omega }T)^n(\eta )\eta ^{-1} =  (I_{\omega }T)^n(\eta )$. Now, by replacing the map $T$ (i.e. the $t$-conjugation) with $I_{\omega }T$, and denoting the new map still with $T$, we can assume that all the elements $T^n(\eta ), n\in \Z $ commute. Then we obtain a contradiction to the bi-orderability as in the Generic Case. Namely, let again $H$ be the subgroup of $G_1$ generated by the subset  $\{T^{n}\eta  \ | \ n\in \mathbb{Z}\}$. It is not necessarily true any more that all these elements are comparable but $H$ is Abelian now. Let again $L = H\cap G_2$. Then we can apply Lemma \ref{thm:use} by taking $V = H, V_0 = L$ and $x = \eta $ to obtain a contradiction.
  
 \medskip
 
 \begin{rem} Let us point out that treating the generic case separately was not necessary. We nevertheless wanted to point out that unless we are unlucky (i.e. if we are in non-generic case), there is no need for the notion of a completion of a bi-ordered group.
 \end{rem}
 
 \section{Non-bi-orderability of $\pi _1(7_6)$}

 In this section we will compute a presentation for the knot group of $7_6$. In Figure 2, we have shown a picture of the Seifert surface $\Sigma $ of the knot $7_6$. Let $K$ denotes the boundary of this surface. In Figure 3, choosing a base point for $\pi _1(\mathbb{S}^3\backslash K)$ and making some convenient choices for the generators of $\pi _1(\mathbb{S}^3\backslash K)$ and $\pi _1(\Sigma \times \{1\})$, we obtain a presentation for the fundamental group $\pi _1(\mathbb{S}^3\backslash K)$. In the figure, the loops $\alpha , \beta , \gamma , \delta $ generate $\pi _1(\Sigma \times \{1\})$, and the loops $a, b, c, d$ generate $\pi _1(\mathbb{S}^3\backslash K)$.
 
 \medskip
 
 We obtain a presentation $$\langle t, \alpha , \beta , \gamma , \delta \ | \ t\alpha t^{-1} = \alpha _{-}, t\beta t^{-1} = \beta _{-}, t\gamma t^{-1} = \gamma _{-}, t\delta t^{-1} = \delta _{-}\rangle .$$ 
 
 Here, the loops $\alpha , \beta , \gamma , \delta $ generate $\pi _1(\Sigma \times \{1\})$, and the loops $a, b, c, d$ generate $\pi _1(\mathbb{S}^3\backslash K)$; $\alpha _{-}, \beta _{-}, \gamma _{-}, \delta _{-}$ are the loops on $\Sigma \times \{-1\}$ parallel to $\alpha , \beta , \gamma , \delta $ respectively. One can compute these loops and obtain a presentation $$\langle t, a, b, c, d \ | \ tat^{-1} = ab, tb^{-1}ac^{-1}t^{-1} = b^{-1}, tct^{-1} = b^{-1}d^{-1}, tdt^{-1} = cd\rangle $$ for the knot group $\Gamma $ of $7_6$. Thus, again $\Gamma $ is isomorphic to a semi-direct product $\mathbb{Z}\ltimes \mathbb{F}^4$ where the acting group $\mathbb{Z}$ is generated by $t$, and the free group $\mathbb{F}_4$ is generated by $a, b, c, d$. 
 
 \medskip

 \begin{figure}[h!]
  \includegraphics[width=6in,height=5in]{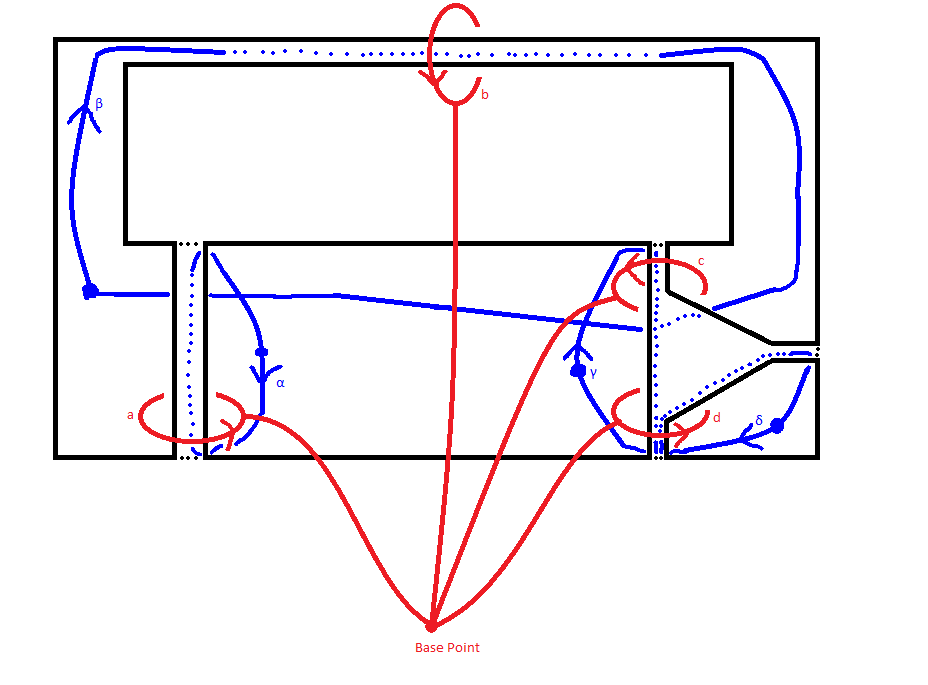}
\caption{A Seifert surface of $7_6$. Blue loops represent the generators of the fundamental group of the surface, while the red loops represent the generators of the fundamental group of complement of the surface.}
\end{figure}

 \medskip
 
 Interestingly, the case of the knot $7_6$ is algebraically identical to case of the knot $6_2$. Indeed, the presentation of $\Gamma $ can also be written as $$\langle t, a, b, c, d \ | \ tat^{-1} = ab, tbt^{-1} = abdb^2, tct^{-1} = b^{-1}d^{-1}, tdt^{-1} = cd\rangle ,$$ from where it is easy to compute that the quotient $\Gamma /[G,G]$ is isomorphic to $\mathbb{Z}\ltimes _{M}\mathbb{Z}^4$ where the acting automorphism is given by the matrix \begin{displaymath}  M =  \left ( \begin{array}{cccc} 1 &  1 & 0 & 0 \\  1 & 3 & -1 & 0 \\  0 & 0 & 0 & 1 \\ 0 & 1 & -1 & 1  \end{array}   \right ). \end{displaymath}
 
 By computing the characteristic polynomial of this matrix we obtain $f(\lambda ) = \lambda ^4 - 5\lambda ^3 + 7\lambda ^2 - 5\lambda +1$.  This polynomial\footnote{Again, it is not accidental that this polynomial coincides with the Alexander polynomial of $7_6$} has two positive irrational real roots of multiplicity 1, and two imaginary irrational roots. By Remark \ref{thm:biorderabilityremark} the group $\mathbb{Z}\ltimes G/[G,G]$ is indeed bi-orderable; thus, as in the case of $6_2$, we do not obtain a contradiction at this level. However, considering the derived series $G_n, n = 0, 1, 2, \dots $ of $G$, we will look at the action of $\mathbb{Z}$ on $G_1/G_2$. To compute this action, let $$[a,b] = u, [a,c] = v, [a,d] = w, [b,c] = x, [b,d] = y, [c,d] = z.$$
 
 Let also $H$ be the subgroup of $G_1$ generated by these six elements. Thus the quotient $H/[H,H]$ is isomorphic to $\mathbb{Z}^6$ generated by $u,v,w,x,y,z$. Then we have $$tut^{-1} = u^2wy, \ tvt^{-1} = u^{-1}w^{-1}y^{-1}, \ twt^{-1} = vwxy,$$ \  $$txt^{-1} = u^{-1}w^{-1}y^{-2}, \ tyt^{-1} = vwx^3y^3z^{-1} , \ tzt^{-1}= x^{-1}y^{-1}z.$$
 
 Thus the action is given by an integral matrix \begin{displaymath}  N =  \left ( \begin{array}{cccccc} 2 & 0 & 1 & 0 & 1 & 0 \\  -1 & 0 & -1 & 0 & -1 & 0 \\  0 & 1 & 1 & 1 & 1 & 0 \\ -1 & 0 & -1 & 0 & -2 & 0 \\ 0 & 1 & 1 & 3 & 3 & -1 \\ 0 & 0 & 0 & -1 & -1 & 1  \end{array}   \right ). \end{displaymath} 
 
 The characteristic polynomial of this matrix has a double root equal to 1, and 4 other imaginary roots.  The 4-dimensional eigenspace corresponding to the 4 imaginary roots contains a non-zero rational vector. Thus the situation is identical to the case of the knot $6_2$ and we conclude that $7_6$ is also non-bi-orderable.

 \bigskip
 
 {\em Acknowledgment:} We are grateful to Dale Rolfsen and Adam Clay for their very close interest and for bringing up several references to our attention.

 \end{document}